\documentclass[12pt]{amsart}
\usepackage{amssymb}

\usepackage{latexsym,amsxtra,amscd,ifthen}
\usepackage{amsmath}
\usepackage{amssymb}

\def\la{\langle}
\def\ra{\rangle}

\def\dim{\mbox{dim\,}}
\def\e{\varepsilon}
\def\epr{$\Box$}           
\def\ext{\mbox{Ext\,}}

\def\gr{\mbox{\bf gr\,}}

\def\Q{{\mathbb{Q}}}

\def\tor{\mbox{Tor\,}}
\def\Z2{\mathbb{Z}_{2}}
\def\Zpl{\mathbb{Z}_{+}}

\newtheorem{definitia}{Definition}
\newtheorem{prop}{Proposition}[section]
\newtheorem{theorema}[prop]{Theorem}
\newtheorem{theorem}[prop]{Theorem}
\newtheorem{lemma}[prop]{Lemma}

\newtheorem{cor}[prop]{Corollary}

\def\proof{\noindent {\bf{Proof} }}

\begin{document}

\author{Dmitri Piontkovski}%
\address{Central Institute of Economics and Mathematics\\
Nakhimovsky prosp. 47, Moscow 117418,  Russia}
\email{piont@mccme.ru}

\author {Sergei D. Silvestrov}
\address{ Centre for Mathematical Sciences, Department of
Mathematics,
Lund Institute of Technology, Lund University,
Box 118, SE-221 00 Lund, Sweden }
\email{sergei.silvestrov@math.lth.se}

\date{August 18, 2009}

\title{Cohomology of 3-dimensional color Lie algebras}

\begin{abstract}
We develop the cohomology theory of color Lie superalgebras
due to Scheunert--Zhang in a framework of nonhomogeneous quadratic Koszul
algebras. In this approach, the Chevalley--Eilenberg complex of a color Lie algebra
becomes a standard Koszul complex for its universal enveloping algebra.
As an application, we calculate cohomologies with trivial coefficients
of $\Z2^n$--graded 3--dimensional color Lie superalgebras.
\end{abstract}

\subjclass[2000]{17B75, 17B56, 16S37}

\maketitle




\section*{Introduction}

\label{sec:intro}

The main objects studied in this paper are the unital associative algebras which are
universal enveloping algebras of color Lie algebras also called generalized Lie algebras, or,
more precisely, $\Gamma$-graded $\varepsilon$-Lie algebras, or generalized Lie algebras.
Since the 1970's, generalized (color)
Lie algebras have been an
object of constant interest in both
mathematics  and  physics
\cite{Agrav,BahtMPZ,cpo,CNS,Gr-Jar},
\cite{Kac-art-adv}--\cite{P}, \cite{R-W(GL)},
\cite{S}--\cite{S-clas}.
The main goal  of the paper is to describe the cohomology with scalar coefficients
of three-dimensional (three generators) Lie color algebras.
Our attention is concentrated on those three-dimensional color Lie algebras
which are in an essential way graded as Lie color algebras by groups larger than
$\mathbb{Z}_{2}$.
The classification of such color Lie algebras has been given in \cite{S-PhD,S-clas}
and presented in terms of commutation relations between
generators. We use this classification in this paper for computing the cohomology.
We believe our results to be useful for better understanding of the structure of the variety of color Lie algebras and their representation theory, and also for shedding further light on
a general constructive approach to cohomology computation.
This work expands the previous works
on three-dimensional color Lie algebras and their representations,
\cite{LS}, \cite{O-Sil}, \cite{PSS}, \cite{SilSig}, \cite{S-stm}
where quadratic central elements and
involutions (real forms) for these algebras have been computed, and
representations for the graded analogues of the Lie
algebra $\mathfrak{sl}(2\, ;\,\mathbb{C})$, of the Lie
algebra of the group of plane motions, and of the Heisenberg Lie algebras, the  three of the non-trivial algebras from
the classification, have been considered.
It is expected that the representation theory for the three-dimensional color Lie algebras,
when developed for all of them,
will show deep connections
with the cohomology structure for universal enveloping algebras of these algebras described
in the present work.

The (co)homology theory for color Lie algebras has been introduced
in~\cite{ShZh-CohomGLA} (see also~\cite{Sh-CohomGLA}).
The theory becomes analogous to the cohomology theory of the usual and $\mathbb{Z}_{2}$--graded Lie algebras;
namely, there are cochain complexes, long exact sequences, homological description of deformations and universal central extensions etc.
In this paper we envelope this cohomology theory in a framework of Koszul algebras. The universal enveloping algebras of color Lie superalgebras are Koszul, and the   Chevalley--Eilenberg complex
becomes a standard Koszul complex for these associative algebras. This gives a kind of Koszul duality for different types of finitely  generated color Lie superalgebras (a type is uniquely defined by the grading and a commutator factor), and
a direct way to calculate cohomology.
We use these techniques to calculate the cohomology (with trivial coefficients)
for 3--dimensional color Lie algebras as classified in~\cite{S-clas}.

The paper is organized as follows. In section~\ref{one}, we give the definition of color Lie algebras and introduce basic notations. In the next section, we first discuss their (co)homologies. Notice that we use another definition of cohomology than the one given in~\cite{ShZh-CohomGLA}, that is, we calculate it via the $\ext$'s of the universal enveloping algebras. However, both these definitions are equivalent (at least over the fields of zero characteristics), according to~\cite[Remark~3.1]{Sh-CohomGLA}
(for a detailed proof, see \cite[Theorem~4]{cpo}).
In the next subsection~\ref{subs-kos}, we briefly recall the theory of Koszul algebras and show
that the universal enveloping algebras of the color Lie algebras are (nonhomogeneous) Koszul.
In section~\ref{abel}, we study the case of abelian color Lie algebras.
Their universal enveloping algebras are homogeneous Koszul, and the classical Koszul duality
gives the duality on the set of their isomorphism classes.
This gives us an easy way to calculate the cohomology of abelian color Lie algebras with scalar coefficients
and to establish some properties of their homologies with coefficients in cyclic modules.
In section~\ref{nonabel}, we use this to describe the method to calculate the algebra of cohomology of arbitrary finite--dimensional  color Lie algebra as a homology of a differential graded algebra, that is, the Koszul dual to the universal enveloping
algebra of its associated abelian color Lie algebra with the differential dual to the multiplication.
In section~\ref{sect-last},
we apply this method to calculate the cohomologies with scalar coefficient of three--dimensional color Lie algebras.

\subsubsection*{Acknowledgment}

The first author was partially supported by Russian Foundation of Basis Research project 02-01-00468.
The research was also partially supported by the Crafoord foundation and The Royal Swedish Academy of Sciences.
The authors gratefully acknowledge
hospitality of Mittag--Leffler Institut during
the Non-commutative Geometry Program 2003-2004, where the main part of this research was performed. The first author is  grateful to the Department of Mathematics, Lund Institute of Technology for hospitality during his visit their.
The authors thank also Daniel Larsson for many useful comments.

The authors specially thank Pasha Zusmanovich  who have found an
error in the resulting table of the previous version of this
paper.

\section{Basic notions}

\label{one}
Given a commutative group $\Gamma$ which will be in what follows referred to as the grading group,
a commutation factor on $\Gamma$ with values in the multiplicative group $k\setminus \{0\}$ of a field $k$ of
characteristic $0$ is a map $$\varepsilon : \Gamma \times \Gamma \mapsto k\setminus \{0\}$$ satisfying three properties:
\begin{eqnarray*}
\varepsilon(\alpha + \beta, \gamma) = \varepsilon(\alpha,\gamma) \varepsilon(\beta,\gamma), \\
\varepsilon(\alpha,\beta + \gamma) = \varepsilon(\alpha,\beta) \varepsilon(\alpha,\gamma), \\
\varepsilon(\alpha,\beta) \varepsilon(\beta,\alpha)=1.
\end{eqnarray*}
The first two properties simply say that the map is a bihomomorphism or a bilinear map on $\Gamma$.

A \mbox{$\Gamma$-graded} $\varepsilon$-Lie algebra
(or a color Lie algebra) is a $\Gamma$-graded linear
space $$ X=\bigoplus_{\gamma \in \Gamma}
{X_{\gamma}} $$ with a bilinear multiplication
(bracket)
$\langle\,\cdot \, , \cdot\,\rangle: X \times X \to X$
obeying

\vbox{\begin{description}
 \item[Grading axiom]
$\langle\, X_{\alpha}\, , X_{\beta}\,\rangle \subseteq  X_{\alpha + \beta}$
 \item[Graded skew-symmetry]
$\langle\, a\, , b\,\rangle=- \varepsilon(\alpha,\beta)\langle \,b\, , a\,\rangle$
 \item[Generalized Jacobi identity]
\begin{equation*}
\hspace{-1.3cm} \varepsilon(\gamma,\alpha)\langle\, a\, , \langle\, b\, , c\,\rangle\,\rangle \mbox{} +
\varepsilon(\beta,\gamma) \langle\, c\, , \langle\, a\, , b\,\rangle\,\rangle
\mbox{} + \varepsilon(\alpha,\beta) \langle\, b\, , \langle\, c\, , a\,\rangle\,\rangle=0
\end{equation*}
\end{description}}
\noindent for all
$\alpha,\beta,\gamma \in \Gamma$  and
\mbox{$a \in X_{\alpha}$ , $b \in
X_{\beta}$ , $c \in X_{\gamma}$}.
The elements of
$\bigcup_{\gamma \in \Gamma} {X_{\gamma}}$ are called
homogeneous.

Any color
Lie algebra $X$ can be embedded in its universal
enveloping algebra $U(X)$ in such a way that,
for homogeneous \mbox{$a \in X_{\alpha}$} and
$b \in X_{\beta}$, the bracket $\langle\,\cdot \, , \cdot\,\rangle$
becomes  the $\varepsilon$-commutator
$\langle\,a, b\rangle\,=
ab-\varepsilon(\alpha,\beta)ba$ which in particular
becomes commutator $[a\, ,b]=ab-ba$
when $\varepsilon(\alpha,\beta)=1$ and
anticommutator $\{a\, ,b\}=ab+ba$ when
$\varepsilon(\alpha,\beta)=-1$ (see \cite{Sh-GLA}).

Before proceeding with the main results and constructions of this work, let us
consider an example of one of the non-commutative algebras appearing in the classification, the color analogue of
the Heisenberg Lie algebra, which is a three-dimensional
\mbox{$\mathbb{Z}_{2}^{3}$-graded} generalized Lie algebra.
Its universal enveloping algebra
is the algebra with three generators $e_1$, $e_2$ and $e_3$ satisfying
defining commutation relations
\begin{align}
e_1 e_2+e_2 e_1&=e_3,  \notag \\
e_1 e_3+e_3 e_1&=0,  \label{rel:colLHeis} \\
e_2 e_3+e_3 e_2&=0       \notag.
\end{align}
When  anticommutators in the left-hand side of (\ref{rel:colLHeis})
are changed  into commutators, we indeed have the
relations between generators in the universal enveloping algebra of
the Heisenberg Lie algebra.

This algebra can be seen as the universal enveloping algebra of a $\mathbb{Z}_{3}^{n}$-graded $\varepsilon$-Lie algebra
where the commutation factor is given by
$\varepsilon(\alpha,\beta) = (-1)^{\alpha_1\beta_1 + \alpha_2\beta_2 + \alpha_3\beta_3}$,
for all  $\alpha=(\alpha_{1},\alpha_{2},\alpha_{3})$  and $\beta=(\beta_{1},\beta_{2},\beta_{3})$ in $\mathbb{Z}_{3}^{n}$.
Now take $X$ to be a \mbox{$\mathbb{Z}_{2}^{3}$-graded} linear space
$$X = X_{(1,1,0)} \oplus X_{(1,0,1)}
\oplus X_{(0,1,1)}$$ with the homogeneous basis
 $e_1 \in X_{(1,1,0)}$, $e_2 \in X_{(1,0,1)}$,
$e_3 \in X_{(0,1,1)}$. The homogeneous
components graded by the elements of
$\mathbb{Z}_{2}^{3}$ different from  $(1,1,0)$, $(1,0,1)$ and
$(0,1,1)$  are zero and so are omitted.
If the \mbox{$\mathbb{Z}_{2}^{3}$-graded} bilinear
multiplication $\langle\,\cdot \, , \cdot\,\rangle$ turns $X$ into
a \mbox{$\mathbb{Z}_{2}^{3}$-graded} generalized Lie algebra,
then with the above commutation factor
$ \langle\, e_{i}\, , e_{i}\,\rangle=0 , \ i = 1,2,3$ and
\begin{equation*}
\langle\, e_{1}\, , e_{2}\,\rangle=c_{12} e_3, \quad
\langle\, e_{2}\, , e_{3}\,\rangle=c_{23} e_1, \quad
\langle\, e_{3}\, , e_{1}\,\rangle=c_{31} e_2 \ .
\end{equation*}
When $a$ and $b$ are in different homogeneous subspaces,
it follows that $\langle\, a\, , b\,\rangle=\langle\, b\, , a\,\rangle$,
whereas $\langle\, a\, , b\,\rangle=-\langle\, b\, , a\,\rangle$
if $a$ and $b$ belong to the same one.
Now put $c_{12}=1$, $c_{23}=0$ and $c_{31}=0$.
The $\mathbb{Z}_{3}^{n}$-graded $\varepsilon$-Lie algebra $X$ so defined  is the color analogue of
the Heisenberg Lie algebra having the algebra defined by relations
\eqref{rel:colLHeis} as its universal enveloping algebra.

The color analogue of the  Heisenberg
Lie algebra \eqref{rel:colLHeis} is one of the algebras in the classification
of three-dimensional color Lie algebras with injective commutation factor
the class of algebras
we work with in this paper.
The condition of injectivity for a commutation factor simply means
that for each fixed $\alpha$ corresponding to a non-zero homogeneous component
the maps $\varepsilon(\alpha, \cdot): \Gamma \rightarrow k\setminus \{0\}$ are injective.
In terms of the grading and commutation relations this means that
one can not reduce the number of non-zero homogeneous components
by simply adjoining some of them into one component of higher dimension
without changing the algebra structure.
For example in case of the color analogue of the
Heisenberg Lie algebra we have three one-dimensional homogeneous components.
The commutation factor can be conveniently represented by a $3\times 3$ matrix
$$
\left(
\begin{array}{rrr}
1 &-1 &-1 \\
-1& 1 &-1 \\
-1&-1 & 1
\end{array}
\right)
$$
indicating the values of the commutation factor, that is the type of the
bracket, between the three non-zero homogeneous components, with $1$ corresponding to a commutator
and $-1$ to an anti-commutator. In terms of the matrix of the commutation factor
the injectivity simply means that there are no equal rows (and by symmetry columns)
in the matrix.

The following simple lemma implies, that together with the known cases of
Lie algebras and superalgebras, the color Lie algebras  with injective commutation factor and
three homogeneous components considered in this paper cover all three-dimensional color Lie algebras
with commutation factors
taking values in $\{ 1, -1\}$.

\begin{lemma}
\label{lemma-inj}
Let $g$ be a color Lie algebra having at most two nonzero homogeneous components,
$$
g = g_i \oplus g_j, \quad i,j \in \Gamma.
$$
 Then either $g$ is abelian or $g$ is a usual
 Lie superalgebra.
  \end{lemma}

 \proof
 Obviously, if the whole $g$ is concentrated in the same component $g_i$, then it is either
 abelian or is an ordinary Lie algebra $g = g_0$.

 First, consider the case $i=0$. Then $\varepsilon (j,j) = \varepsilon (0,j) \varepsilon (j,j)$.
 Hence $\varepsilon (0,j) =1$. Because $\varepsilon (0,0) =1$,   we have that $g$ is either
 a Lie algebra or a $\mathbb{Z}_{2}$--graded  Lie algebra depending on the sign of $\varepsilon (j,j) = \pm 1$.

 The case $j =0$ is completely analogous.
 So, we may assume that $i , j\ne 0$.
 Since $i + j \ne i,j$, we have $\la  g_i ,g_j \ra \subseteq g_{i + j} = 0$.
 Suppose that  $g$ is not abelian, that is, say, $\la  g_i ,g_i \ra \ne 0$.
 Then $\la  g_i ,g_i \ra \subseteq g_j$, so, $ j = i+i$. Thus,
 $\varepsilon (i,j) = \varepsilon (i,i)^2 =1$, and hence
 $\varepsilon (j,j) = \varepsilon (i,j)^2 =1$.
 Then  $g$ again becomes $\mathbb{Z}_{2}$--graded, with odd component
 either $g_i$ or zero, depending on the sign of $\varepsilon (i,i) = \pm 1$.
\epr


\section {Koszul algebras and color Lie algebras}

In this section, we describe the general construction and methods of cohomology theory
of color Lie algebras. Here we assume that the ground field $k$
has characteristics different form  $2$ and $3$.

\subsection {(Co)homology of color Lie algebras}

In this subsection, we first consider associative $k$--algebras and modules over them.
All algebras are assumed to be {\it augmented}, that is,
admitting a surjective homomorphism (augmentation) $A \to k$;
its kernel (augmentation ideal) is denoted by $A_+$.
All modules over associative algebras below are assumed to be  right-sided.

Let us introduce some short notations for (co)homologies of modules over augmented associative
algebras. For a module $M$ over a $k$--algebra $R$, we will denote by $H_i (M) = H_i (R,M)$ its
homology, i.~e., the vector space $\tor_i^R (M,k)$. By $H_i A$ we will denote
the vector space $\tor_i^R (k,k)$. In the same way, we define
$H^i (M) = H^i (R,M) := \ext^i_R (k,M)$  and
$H^i R := H^i (R,k) := \ext^i_R (k,k)$.

\begin{definitia}
Let $g$ be a (color) Lie algebra, let $R = U(g)$  be its universal enveloping algebra, and let
$M$ be an $R$--module.
By the $i$--th {\it homology} of $g$ with coefficients in $M$ we will mean the vector space $H_i (g, M)
:=H_i (R,M)$, by $i$--th {\it cohomology} we will mean     $H^i (g, M)
:=H^i (R,M)$. By $i$--th (co)homology of $g$ with trivial coefficients we will mean
the vector spaces $H_i g := H_i (g, k) = H_i R$ and   $H^i g :=  H^i (g, k) = H^i R$.
\end{definitia}

Because of~\cite[Remark~3.1]{Sh-CohomGLA} (see also~\cite[Theorem~4]{cpo}), these
cohomology groups are the same as the ones defined in~\cite{ShZh-CohomGLA}, at least in the case
when the ground field is of zero characteristics.

\subsection {Koszulity}
\label{subs-kos}

We will call a linear space over $k$, a $k$--algebra, or
$k$--algebra  module {\it positively graded,} if it is ${\bf Z}_+$--graded
and finite--dimensional in every component.
A positively graded associative algebra $R= R_0 \oplus R_1 \oplus \dots$
is called {\it connected,} if its zero component $R_0$
is $k$; a connected algebra is called {\it standard,}
if it is generated by $R_1$ and a unit.
All modules over positively graded algebras below are assumed
to be positively graded and right-sided.
For a positively graded module $M$ over a connected algebra, its (co)homology vector space is again graded,
that is, $H_i(M) = \oplus_{i}H_i(M)_j $,  $H^i(M) = \oplus_{i}H^i(M)^j $
 and for (co)homologies with trivial
coefficients we have the following duality isomorphisms of graded vector spaces:
$$ H_i(A) = H^i(A)^*.$$

\begin{definitia}
A graded module $M$ over a standard associative algebra is called

--- linear (of degree $d$), if it is generated in degree $d$, i.e.,
    $H_0(M)_j = 0$ for $j \ne d$;


--- quadratic, if it is linear of degree $d$ and all
    its relations may be choosen in degree $d+1$, i.e., $H_1(M)_j = 0$ for $j \ne d+1$;

--- Koszul, if it has linear free resolution, i.e.,
    $H_i(M)_j = 0$ for all $i\ge 0, j \ne i+d$.
\end{definitia}


\begin{definitia}{\cite{pri}}
A standard algebra $R$ is called (homogeneous) Koszul if the trivial module $k_R$ is Koszul, i.e.,
every homology module $H_i R$  is concentrated in degree $i$.
\end{definitia}

Notice that the notion of homogeneous Koszul algebra has a lot of applications and a wide
theory (see \cite{bf, pp} and references therein).

Priddy in his original paper \cite{pri} had also introduced the
 notion of Koszulity for  some
non-graded algebras.
 Let $A$  be an augmented associative
algebra minimally generated by a finite-dimensional
vector space $V$  (that is, the augmentation
ideal is generated by $V$). Then $A$ is called {\it almost quadratic} if its relations are
nonhomogeneous noncommutative polynomials of degree two, that is, $A$  is isomorphic to a quotient
algebra of the free algebra $F = T ( V) $ by an ideal $I$  generated by a subspace of relations
$R \subset F_1 \oplus F_2 = V \oplus V \otimes V$.

 There are two graded algebras
associated to $A$. The first one, denoted here by $A^0$, is a quotient algebra
of $A$ by an ideal generated by the projection of $R$ on $ F_2 = V \otimes V$.
It is homogeneous quadratic.
Also, the algebra $A$ is obviously
filtered by total degree, and the associated graded algebra $\gr A$ is again a standard algebra generated by $V$.

In these terms, the definition of Priddy looks as follows.
\begin{definitia}{\cite{pri}}
An almost quadratic  algebra $A$ is called Koszul if both graded algebras $A^0$ and  $\gr A$
 are isomorphic
to each other and Koszul.
\end{definitia}

A discussion on this concept could be found in~\cite{bg,pp}.

Priddy had also introduced the following special kind of Koszul algebras.
We extend his definition to our nonhomogeneous case.
In modern terms, the definition looks as follows.

\begin{definitia}
An almost quadratic  algebra $A = F/I$ is Poincar\'e--Birkhoff--Witt (PBW for short), if, for some
choice of bases $v_1, \dots, v_n$ of $V$ and $r_1, \dots, r_m$ of $R$, the set $ r = \{ r_1, \dots,
r_m \}$ is a Groebner basis of the ideal $I$ w.~r.~t. the degree--lexicographical order with $v_1 >
\dots > v_n$.
\end{definitia}

\begin{theorem}
Every almost quadratic PBW algebra is Koszul.
\end{theorem}

\proof
By theorem of Priddy \cite{pri}, all graded (=homogeneous) PBW algebras are Koszul.

Let $r_i = f_i + l_i$ for $i = 1, \dots, m$,
where $l_i \in V$ is the linear part and $f_i \in V \otimes V$ is the homogeneous quadratic part,
and let $t_i$ be the leading monomial of $r_i$.
By the Diamond Lemma (see~\cite{ufn}), the  condition that $r$ is a Groebner basis means
that, if $t_i v_p = v_q t_j$ for some $i,j,p,q$,
 then the $s$--polynomial
$s = r_i v_p - v_q r_j$ may be reduced to zero w.~r.~t. $r$.
In particular, its high--degree part (of degree 3) $s' = f_i v_p - v_q f_j$
may be reduced to an element of degree less than 3 w.~r.~t. $r$, that is, it is reducible
to zero w.~r.~t. the set $f = \{ f _1, \dots, f_m \}$ of relations of $A^0$.
This means that $f$ is Groebner basis of relations of $A^0$; in particular, the
previously mentioned result of Priddy implies that $A^0$ is Koszul. Also, by the definition of Groebner bases,
the algebras $A$ and $A^0$ have the same linear basis consisting of
all the monomials on generators without submonomials equal to $t_i$.
Because the degree--lexicographical order is an extension of the
partial order $``<<"$ given by degrees of elements of $F$, the associated graded algebra
$\gr A$ w.~r.~t. the filtration induced by $``<<"$
has the same linear basis. Thus, it is isomorphic to $A^0$.
\epr

 The terminology is due to the fact
that the universal enveloping algebras of finite--dimensional Lie (super)algebras are PBW.
In this connection, the following re-formulation of
Poincar\'e--Birkhoff-Witt theorem~\cite{Sh-GLA,Sh-GTC}
for color
Lie algebras is not surprising.

\begin{theorema}
\label{PBW}
Let $R = U (g)$  be universal enveloping algebra of a finite--dimensional
color Lie algebra $g$. Then the almost quadratic algebra $R$ is PBW (in particular, Koszul).
The associated graded quadratic algebra is the universal enveloping algebra of some
abelian color Lie algebra.
\end{theorema}

\proof
Let  $v_1, \dots, v_n$ be a $\Gamma$--homogeneous basis of $V$, let $r$ be the set
$f_{ij} = v_i v_j + {\e (v_i, v_j) } v_j v_i - \sum_k c_{ij}^k $,
and let $G$ be the set of the leading monomials of these elements.
By the PBW-Theorem~\cite{Sh-GLA,Sh-GTC},
the set $B$ of monomials on the variables     $v_1, \dots, v_n$ whose all degree two
submonomials do not belong to $G$, is linear independent modulo $I$.
By definition, this means that
the set $r$ is a Groebner basis of $I$.
It follows that the associated graded algebra has the relations
 $h_{ij} = v_i v_j + {\e (v_i, v_j) } v_j v_i$. By definition, it is the universal enveloping algebra
 of the  abelian color Lie algebra with the same  $\Gamma$--homogeneous basis
 and the same commutation factor as $g$.
\epr

\begin{cor}
\label{PBW ab}
Let $g$ be an abelian color Lie algebra. Then its universal enveloping algebra  $R = U (g)$ is
homogeneous Koszul.
\end{cor}

\section {The case of abelian color Lie algebras}

\label{abel}

\subsection {Koszul duality}
\label{kosdual}

For any standard quadratic algebra $A$ with the space of generators $V = A_1$ and the
space of relations $R \subset V \otimes V$, its quadratic dual algebra
$A^!$ is, by definition, a standard quadratic algebra generated by the space $V^*$  dual to $V$
with the set of relations $R^\perp \subset V^* \otimes V^*$ (the annihilator of $R$ in the dual space to $V \otimes V$).
The importance of this algebra is
due to the following fact~\cite{pri}: if the algebra $A$ is Koszul, then $H_i A \simeq H^i A \simeq
A^!_i$.

Let $V$ be a $n$--dimensional vector space spanned by $v_1, \dots, v_n$.
Let $[n]$
be the $n$--element set $\{ 1, \dots, n \}$ and $P$ be a set of unordered pairs of integers
$\{  (i,j) | i, j \in [n], i < j\}$. For any two subsets $J \subset [n], Q \subset P$
let us define a standard associative algebra $A_{J,Q} = T \la V\ra / I_{J,Q}$, where the ideal $I_{J,Q}$
is generated by the following set of quadratic elements:
$$
    v_i v_j - v_j v_i, (i,j) \in Q;
$$
$$
    v_i v_j + v_j v_i, (i,j) \in P - Q;
$$
$$
    v_i^2, i \in J.
$$

Such relations can be put into a one-to-one correspondence with
unordered graphs. Let $\Gamma = (V,E)$ be a graph
with a set of vertices $V=\{v_i\}_{i\in [n]}$ and a set of edges $E$, where $v_i$
and $v_j$ are connected by an edge if and only if they
anti-commute, that is if $(i,j) \in P - Q$ or $i=j \in J$
(if the coefficient
field is not of characteristic $2$, then $i=j \in J$ is equivalent
to anti-commutativity of generator with itself).

\begin{prop} The following relation holds:
$$ A_{J,Q}^! =  A_{[n]-J,P-Q}.
$$
\end{prop}

The associative algebra $A_{J,Q} = T ( V) / I_{J,Q}$ can be also
viewed as a universal enveloping algebra of a commutative
$\mathbb{Z}_{2}^{n}$-graded (color) generalized Lie
algebra.

The following statement follows from the constructions in~\cite{S}.

\begin{prop}
     Let $g$ be an $n$--dimensional abelian color Lie algebra with commutation
     factor taking the values in  the set $\{ -1,1\}$. Then its universal enveloping
algebra  $A = U (g)$ is isomorphic to the algebra $A_{J,Q}$ for some $J,Q$.

Conversely, every algebra   $A_{J,Q}$ with $n$ generators is isomorphic to the  universal enveloping
algebra of some $n$--dimensional abelian color Lie algebra.
\end{prop}

\subsection {(Co)homology with trivial coefficients}
\label{nonabel}

In the view of the previous subsection, the dimensions of (co)homologies
of any abelian color Lie algebra $g$ are equal to the
dimensions of homogeneous components of suitable PBW-algebra $A_{J,Q} = (U(g))^!$.
To calculate these dimensions, it remains to show a formula for the generating function for these dimensions, that is, for the Hilbert series
$$ A_{J,Q} (z) = \sum_{i=0}^\infty z^i \dim (A_{J,Q})_i.
$$

\begin{prop}
Let $q$ be the cardinality of elements in the set $Q$.Then
$$ A_{J,Q} (z) = \frac{(1+z)^{q}}{(1-z)^{n-q}}.
$$
\end{prop}

\proof
Since the algebra $A_{J,Q}$ is PBW, as a graded vector space it is isomorphic to the span
of the monomials which has no divisors of the form
$\{ v_i v_j | i < j \} \cup \{ v_i^2 | i \in Q\}$.
It follows that its Hilbert series is equal to the Hilbert series of any algebra spanned by the same
set monomials, in particular, to the one of the algebra $\Lambda (k^q) \otimes S (k^{n-q})$.
Hence
$$
A_{J,Q} (z) = \Lambda (k^q)(z)  S (k^{n-q}) (z) =  \frac{(1+z)^{q}}{(1-z)^{n-q}}.
$$
\epr

\begin{cor}
\label{poincare}
Let $g$ be an $n$--dimensional abelian color Lie algebra such that in its $\Gamma$--homogeneous
relations there are exactly $q$ squares of the generators (that is, there are exactly $q$
relations of the form $\la v_i, v_i\ra = 0$ with $\e (i,i) = -1$). Then its universal enveloping
algebra  $A = U(g)$ has Poincare series of the
following form:
$$
   P_A (z)  :=  \sum_{i=0}^\infty z^i \dim H_i A  =  \frac{(1+z)^{n-q}}{(1-z)^{q}}.
$$
In particular, the homology algebra is finite-dimensional iff there are no squares among the $\Gamma$--homogeneous
relations.
\end{cor}

\subsection {Further Koszul properties}
In this paragraph, we describe some additional properties of the algebras  $A_{J,Q}$.
The subsequent does not depend on them.

A homogeneous Koszul algebra may also have other homological properties
depending on the structure of its ideals and cyclic Koszul modules. Such a theory is connected to the
concept of the so-called {\it Koszul filtration}. For commutative Koszul algebras, this concept has
been studied in several papers~\cite{crv, ctv,i-kos}. The non-commutative version of
the theory has been developed in~\cite{pi}.

\begin{definitia}{\cite{pi}}
Let $R$ be a standard algebra.
A set ${\bf F}$ of degree-one generated right-sided ideals in $R$ is called a Koszul filtration
if $0 \in {\bf F},  R_+ \in {\bf F}$,
and for every $0 \ne I \in {\bf F}$
there are $I \ne J \in {\bf F}$ and $x \in R_1$
such that $I = J + x R$ and the ideal $(x:J) :=
\{ a \in R | xa \in J \}$ lies
in ${\bf F}$.
\end{definitia}
As in the commutative case, every algebra admitting Koszul
filtration is Koszul, as well as every ideal $I \in F$ and its cyclic module $R/I$.

The smallest possible Koszul filtration consisting of the complete flag of  degree-one generated right-sided ideals
$0 = I_0 \subset \dots \subset I_n = R_+$ is called the Groebner flag: it corresponds to the sequence
$x_1, \dots, x_n$ of generators of $V = A_1$.
An algebra having a Groebner flag is called
initially Koszul.

Applying a criterion for initially Koszul algebras given in~\cite{pi}, we obtain the following
results for color Lie algebras.
\begin{prop}
 Let $g$ be a $n$--dimensional abelian color Lie algebra. Then its universal enveloping algebra $A = U(g)$ is initially Koszul.

Moreover, it is initially Koszul w.~r.~t.  any sequence $x_1, \dots, x_n$ of homogeneous generators
of $g$, that is, $A$ is universally initially Koszul~\cite{i-kos} (in particular, it is strongly
Koszul~\cite{hhr}).
\end{prop}

\begin{cor}
Let $W$  be a subset of the set of homogeneous generators of an $n$--dimensional abelian color Lie
algebra $g$, and let $M = U(g) / W U(g)$. Then the graded vector space
$H_i(g,M) $
is concentrated in degree $i$ for all $i \ge 0$.
\end{cor}

\section {The case of arbitrary finite-dimensional color Lie algebras}

\subsection{Priddy's method to calculate the cohomology}

Let $g$ be an $n$--dimensional color Lie algebra, and let $g_{Ab}$ be the associated abelian
color Lie algebra. Denote $A = U(g), \bar A = U(g_{Ab})$ and $ A^! = \bar A^! $.
Let $V = \bar A_1$ be the
set of generators of $A$  which is identified with $g$. The multiplication $\mu : V \otimes V \to V$
in $g$ induces the dual map $d_1: V^* \to  V^* \otimes  V^*$, whose image may be identified with
$A^!_2$. Let us expand $d_1$ up to the differential $d$ of the whole algebra $A^!$  by the Leibniz
rule
$$
     d (ab) = d(a) b + (-1)^{|a|} a d (b),
$$
where $|\cdot |$ denotes
the ${\bf Z}$--grading degree in $A^!$.
Setting $d 1 = 0$, we may consider the algebra $A^!$ as a differential graded algebra
$(A^!, d)$.

The following theorem was proved by Priddy for the case of arbitrary Koszul algebra $A$.

\begin{theorem}
\label{prid}
In the notation above,
$H^\bullet g = H^\bullet (A^!, d) .$
\end{theorem}

The first corollary is well-known for classical Lie algebras.

\begin{cor}[\cite{ShZh-CohomGLA}]
\label{cor_h1}
 Let $g$ be a finite--dimensional color Lie algebra. Then
$$
        H^0 g = k,
$$
$$
         H^1 g  = \left( g/[g,g] \right)^*.
$$
\end{cor}

\section{Main result}
\label{sect-last}

\begin{theorem}
\label{main}
The cohomologies of three--dimensional color Lie algebras over $\mathbb{C}$
with  injective
commutation factor taking values in $\{ -1, 1 \}$
are described in Table~\ref{answer}.
\end{theorem}

{\bf Comments on Table~\ref{answer}}

The number in brackets in the first column denotes the number in the  classification table in~\cite{PSS}. 
Among the 27 cases in~\cite{PSS}, 12 color Lie algebras are abelian. Their homology is described in Corollary~\ref{poincare}. The rest 15 cases are listed in our Table~\ref{answer}.
The graphs in the second column denote the associated graded abelian algebras
as described in subsection~\ref{kosdual}. We do not repeat the graphs in the subsequent lines if the
algebras have the same associated abelian algebras.
The Betti numbers $h_n = \dim H^n (g)$ are given in the last column in the form of their generating function (Poincare series) $P_g(z) = \sum_{n \ge 0} h_n z^n$.
For infinite--dimensional algebras $H^\bullet (g)$, we also give the explicit values for
$h_n , n \ge 1$.

The {\bf proof} of Theorem~\ref{main} is given in Table~\ref{calc_process}.
We calculate the cohomology using Theorem~\ref{prid}.
The numbering of cases is the same as in the previous Table~\ref{answer}.
The linear bases of the differential graded algebras (that is, the universal enveloping algebras
of the dual abelian Lie algebras) for the cohomology calculation are given in the second column.
If an entry in the second column is empty, then it means that it is the same as the preceding first nonempty  entry.
The values of the differentials on the generators $f_i = 2 e_i^* $ (where $\{ e_i^* \}$ is the basis of $g^*$
dual to a fixed basis  $\{ e_i \}$ of $g$) are given in the third column.
In the last three columns, we give the bases of the vector spaces of cocycles $Z^n$, coboundaries $B^n$,
and cohomology $H^n$. Sometimes we omit the case $n =1$,
because the groups $H^1g$ are described in Corollary~\ref{cor_h1}.
Notice that the multiplication structure of the algebra $H^\bullet (g)$ is induced by the multiplication in the
initial differential graded algebra $A^!$.
\newpage


{\tiny
\begin{table}[hp]
   \begin{center}
     \caption{\protect Cohomologies with trivial coefficients\label{answer}}

   \end{center}
\end{table}



\begin{table}[hp]
   \begin{center}
     \caption{\protect Cohomology computation \label{calc_process}}
    %
                                       $\\
\hline
     \end{tabular}
   \end{center}
 \end{table}
}

\end{document}